%% file: main.tex
\documentclass[a4paper, 11pt]{article}
\usepackage{amsmath,amssymb,amsthm,mathtools}
\usepackage{booktabs,enumitem}
\usepackage{graphicx}
\usepackage{xcolor}
\usepackage[colorlinks=true,linkcolor=blue,citecolor=blue,urlcolor=blue]{hyperref}
\usepackage[capitalise]{cleveref}

\newtheorem{theorem}{Theorem}[section]
\newtheorem{lemma}[theorem]{Lemma}
\newtheorem{proposition}[theorem]{Proposition}
\newtheorem{corollary}[theorem]{Corollary}
\theoremstyle{definition}
\newtheorem{definition}[theorem]{Definition}
\newtheorem{construction}[theorem]{Construction}
\newtheorem{example}[theorem]{Example}

\theoremstyle{remark}

\newcommand{\F}{\mathbb{F}}
\newcommand{\C}{\mathbb{C}}
\newcommand{\one}{\mathbf{1}}
\newcommand{\GH}{\mathrm{GH}}

\newcommand{\SRG}{\mathrm{SRG}}
\newcommand{\DDG}{\mathrm{DDG}}
\newcommand{\AG}{\mathrm{AG}}

\newcommand{\tr}{\operatorname{tr}}

\title{Auxiliary matrices and decompositions of affine divisible design graphs}

\author{Keita Ishizuka\thanks{Information Technology R\&D Center, Mitsubishi Electric Corporation, Kanagawa, Japan. \texttt{keitaishizuka1994@gmail.com}}
~and
Hadi Kharaghani\thanks{Department of Mathematics and Computer Science, University of Lethbridge, Lethbridge, Alberta, T1K 3M4, Canada. \texttt{kharaghani@uleth.ca}}}
\date{\today}

\begin{document}
\maketitle

\begin{abstract}
Kabanov constructed divisible design graphs by gluing copies of affine designs along bijections between their parallel classes, arranged according to a symmetric design. On the other hand, auxiliary matrices from finite fields have been used to decompose the complete graph into divisible design graphs and strongly regular graphs of the same kind. We show that Kabanov's graphs are exactly the graphs obtained from auxiliary matrices of affine designs, and we determine their eigenvalues with multiplicities. As in the finite field case, when the symmetric design is trivial, shifting the bijections decomposes the complete graph into Kabanov's graphs and cliques. A union of graphs of such a decomposition is a divisible design graph exactly when the shifts form a difference set, which gives new infinite families of divisible design graphs.


\end{abstract}

\section{Introduction}

Strongly regular graphs are among the central objects of algebraic combinatorics, linking graph theory with design theory, coding theory and finite geometry~\cite{BvM2022}. A \emph{divisible design graph} (DDG) with parameters \((v,k,\lambda_1,\lambda_2,m,n)\) is a \(k\)-regular graph on \(v=mn\) vertices whose vertex set is partitioned into \(m\) classes of size \(n\) such that two distinct vertices have \(\lambda_1\) common neighbours if they lie in the same class and \(\lambda_2\) common neighbours otherwise. If \(\lambda_1=\lambda_2=\lambda\), the graph is strongly regular with parameters \((v,k,\lambda,\lambda)\). Divisible design graphs were introduced in~\cite{HKM2011} as the graphs whose adjacency matrix is the incidence matrix of a symmetric divisible design. They generalize the strongly regular graphs with \(\lambda=\mu\) and, like these, connect design theory with algebraic graph theory; a proper divisible design graph other than \(mK_n\) and its complement is a Deza graph in the sense of~\cite{EFHHH1999}; see~\cite{HKM2011}. Many constructions are known~\cite{HKM2011,PS2022}, and all proper connected divisible design graphs with at most \(39\) vertices have been determined, except for three parameter sets~\cite{PS2022}.

One of the main sources of such graphs is the gluing of affine designs. Wallis~\cite{Wallis1971} constructed strongly regular graphs from affine designs and a Steiner \(2\)-design. Fon-Der-Flaass~\cite{FDF2002} modified the case in which the Steiner \(2\)-design has blocks of size \(2\) and obtained hyperexponentially many strongly regular graphs with the same parameters. Muzychuk~\cite{Muzychuk2007} extended these ideas to all cases of Wallis's construction and showed that the Steiner \(2\)-design can be replaced by a partial linear space, which gave new prolific constructions. Kabanov carried the method over to divisible design graphs. In~\cite{Kabanov2022} he modified the Wallis--Fon-Der-Flaass construction to produce divisible design graphs, in~\cite{Kabanov2026DM} he gave two prolific constructions of infinite series of them, and in~\cite{Kabanov2026JCD} he arranged copies of a collection of affine designs according to a symmetric design: two vertices in different copies are adjacent if one lies in the image, under a prescribed bijection between parallel classes, of the block through the other. Since the bijections can be chosen freely, these constructions produce many graphs with the same parameters. The eigenvalues of these graphs follow from their parameters, but Kabanov remarks that their multiplicities cannot be determined in this way~\cite[after Theorem~3]{Kabanov2026JCD}.

Decomposing the complete graph into strongly regular graphs and divisible design graphs has been a recurring theme in the work of the second author and his coauthors~\cite{KT2003,KSS2017,KharaghaniSuda2018,CIKSS2026}, and the main tool of~\cite{KSS2017,KharaghaniSuda2018,CIKSS2026} has been \emph{auxiliary matrices}: \((0,1)\)-matrices obtained from a generalized Hadamard matrix by replacing each entry by a permutation matrix, an idea going back to the second author's construction of Bush-type Hadamard matrices~\cite{Kharaghani1985}. Kharaghani, Sasani and Suda~\cite{KSS2017} arranged the auxiliary matrices from the multiplication table of \(\F_{2^m}\) by a symmetric Latin square with constant diagonal and decomposed the complete graph on \(2^{2m}(2^m+2)\) vertices into \(2^m\) strongly regular graphs and \(2^m+2\) cliques. For an odd prime power \(q\), Kharaghani and Suda~\cite[Section~5]{KharaghaniSuda2018} obtained in the same way divisible design graphs \(N_\alpha\), \(\alpha\in\F_q\), whose products \(N_\alpha N_\beta\) depend only on \(\alpha-\beta\), and determined the algebra they generate. In~\cite[Appendix~A]{CIKSS2026} shifted auxiliary matrices over \(\F_q\), \(q\) a prime power, were used to decompose the complete graph and to construct Deza graphs. In~\cite{KharaghaniSuda2019} auxiliary matrices are defined for any affine resolvable design. In all these decompositions the auxiliary matrices come from a finite field.

In this paper, we first show that Kabanov's graphs are exactly the graphs given by auxiliary matrices of affine designs (\cref{prop:kabanov}). This gives a matrix proof of their parameters (\cref{thm:ddg}) and determines their eigenvalues with multiplicities (\cref{cor:spectrum}), which were left open in~\cite{Kabanov2026JCD}. Going beyond Kabanov's graphs themselves, we shift the bijections by an abelian group \(G\) of order \(q\) and show that the shifted graphs partition the edges of \(N\otimes J\); for \(N=J-I\) this decomposes the complete graph into Kabanov's graphs and cliques (\cref{thm:decomposition}). More generally, we determine the eigenvalues, with multiplicities, of every union of the shifted graphs (\cref{cor:spectrum}). Finally, a union over \(S\subseteq G\) is a divisible design graph if and only if \(S\) is a difference set (\cref{prop:diffset}). Known difference sets give infinite families of new divisible design graphs; \cref{tab:ddg} lists those with at most \(10^5\) vertices (\cref{cor:families}).

This paper is organized as follows. In \cref{sec:prelim} we fix notation, recall affine designs in matrix form, review the decomposition of the complete graph by auxiliary matrices on the smallest example, and recall Kabanov's construction. In \cref{sec:aux} we identify the auxiliary matrices from finite fields with parallel classes of \(\AG(2,p)\) and bijections between them. In \cref{sec:construction} we give the construction, compute its square and show that it yields exactly Kabanov's graphs. In \cref{sec:shift} we shift the bijections, decompose \(N\otimes J\), determine the spectra, and characterize the unions that are divisible design graphs by difference sets.

\section{Preliminaries}\label{sec:prelim}

Throughout, \(I_n\), \(J_n\) and \(\one_n\) denote the identity matrix, the all-ones matrix and the all-ones column vector of size \(n\); subscripts are omitted when the size is clear. The Kronecker product is \(\otimes\). Rows and columns of a \(q\times q\) matrix are indexed by a set \(Q\) of size \(q\).

\subsection{Permutation matrices}
For a permutation \(\sigma\) of \(Q\) let \(P_\sigma\) be the permutation matrix with \(P_\sigma[c,d]=1\) if and only if \(d=\sigma(c)\). Then
\[
P_\sigma P_\tau=P_{\tau\circ\sigma},\qquad P_\sigma^{T}=P_{\sigma^{-1}},\qquad P_\sigma\one=\one,\qquad P_\sigma J=JP_\sigma=J .
\]
For a prime \(p\) and \(Q=\F_p\), let \(U\) be the circulant matrix with first row \((0,1,0,\dots,0)\) and \(R\) the back-identity matrix, so that
\[
U^h[c,d]=1\iff d=c+h,\qquad (U^hR)[c,d]=1\iff c+d=-1-h .
\]
Thus \(U^h=P_\tau\) for the translation \(\tau(z)=z+h\), and \(U^hR=P_\psi\) for the reflection \(\psi(z)=-1-h-z\).

\subsection{Generalized Hadamard matrices}
A \emph{generalized Hadamard matrix} \(\GH(p,\lambda)\) over \(\F_p\) is a \(p\lambda\times p\lambda\) matrix \(H=(h_{ij})\) with entries in \(\F_p\) such that for any two distinct rows \(i\ne i'\) the differences \(h_{ij}-h_{i'j}\) take every value of \(\F_p\) exactly \(\lambda\) times. The \emph{expansion} of a row \((h_1,\dots,h_N)\) is the \(p\times pN\) matrix \([U^{h_1}\ U^{h_2}\ \cdots\ U^{h_N}]\).

\begin{lemma}\label{lem:ghfacts}
\begin{enumerate}[label=(\alph*)]
\item The matrix \(H_p=(at)_{a,t\in\F_p}\) is a \(\GH(p,1)\).
\item If \(\widehat r\) and \(\widehat r'\) are the expansions of two distinct rows of a \(\GH(p,\lambda)\), then
\begin{equation}\label{eq:gh}
\widehat r\,\widehat r^{T}=p\lambda\, I_p,\qquad \widehat r\,\widehat r'^{T}=\lambda J_p .
\end{equation}
\end{enumerate}
\end{lemma}

\begin{proof}
(a) is standard; see~\cite{deLauney2007}. (b) We have \(\widehat r\,\widehat r^{T}=\sum_jU^{h_j}(U^{h_j})^{T}=p\lambda I\) and \(\widehat r\,\widehat r'^{T}=\sum_jU^{h_j-h'_j}=\lambda(I+U+\dots+U^{p-1})=\lambda J_p\).
\end{proof}

\subsection{Affine designs as matrices}
An affine design with \(n\) points, \(q\) blocks in each parallel class and \(\kappa\) parallel classes is recorded by the incidence matrices of its parallel classes; the incidence matrix of a class is the \(q\times n\) matrix whose rows are indexed by the \(q\) blocks of the class and whose columns are indexed by the \(n\) points. We use the following matrix form; see~\cite{Bose1942} and~\cite{AbelGeYin2007} for affine designs.

\begin{definition}\label{def:system}
Let \(q\ge2\) and \(n\) be integers with \(q^2\mid n\) and \((q-1)\mid(n-1)\), and put \(\kappa=(n-1)/(q-1)\). An \emph{\((n,q)\)-affine system} is a family \(R_1,\dots,R_\kappa\) of \(q\times n\) \((0,1)\)-matrices such that
\begin{enumerate}[label=(A\arabic*)]
\item\label{A1} every column of every \(R_u\) contains exactly one \(1\);
\item\label{A2} \(R_uR_u^{T}=\frac nq I_q\) for every \(u\);
\item\label{A3} \(R_uR_v^{T}=\frac{n}{q^2}J_q\) for \(u\ne v\).
\end{enumerate}
\end{definition}

The rows of \(R_u\) are the blocks of the \(u\)-th parallel class: \ref{A1} says that each class partitions the points, \ref{A2} that the blocks have size \(n/q\), and \ref{A3} that two blocks from different classes meet in \(n/q^2\) points. By \cref{lem:system} below, an \((n,q)\)-affine system is a resolvable \(2\)-design in which any two non-parallel blocks meet in a constant number of points, that is, an affine design. Conversely, the parallel classes of an affine design with \(n\) points and \(q\) blocks in each class satisfy \ref{A1}--\ref{A3} by~\cite[Theorem~5.21]{Stinson2004}, and there are \(r=k+\lambda=(n-1)/(q-1)\) of them by Bose's inequality~\cite{Bose1942}; see~\cite[Theorem~5.14 and Definition~5.17]{Stinson2004}. So an \((n,q)\)-affine system is the same thing as an affine design. Examples are the point-hyperplane designs of \(\AG(m,q)\) (\(n=q^m\)), affine planes of order \(q\) (\(n=q^2\)) and Hadamard \(3\)-designs (\(q=2\)).

\begin{lemma}\label{lem:system}
Let \(R_1,\dots,R_\kappa\) be an \((n,q)\)-affine system and \(\lambda_A=\frac{n-q}{q(q-1)}\). Then \(R_u^{T}\one_q=\one_n\), \(R_u\one_n=\frac nq\one_q\), and
\[
\sum_{u=1}^{\kappa}R_u^{T}R_u=(\kappa-\lambda_A)I_n+\lambda_A J_n .
\]
\end{lemma}

\begin{proof}
The first identity is \ref{A1}, and the second follows from \ref{A2}. If the system comes from an affine design, the matrix \(M=[R_1^{T}\ \cdots\ R_\kappa^{T}]\) is the incidence matrix of a \(2\)-\((n,\frac nq,\lambda_A)\) design with replication number \(\kappa\), and the third identity is \(MM^{T}=(r-\lambda)I+\lambda J\)~\cite[Theorem~1.13]{Stinson2004}. For a general affine system we argue directly. Let \(M=[R_1^{T}\ \cdots\ R_\kappa^{T}]\), so that \(MM^{T}=\sum_uR_u^{T}R_u\). By \ref{A2} and \ref{A3}, \(M^{T}M=\frac nqI+\frac{n}{q^2}(J-I_\kappa\otimes J_q)\), which has the eigenvalue \(\frac nq\) on the vectors whose entries sum to \(0\) on each block of size \(q\), a space of dimension \(\kappa(q-1)=n-1\). Hence the symmetric matrix \(MM^{T}\) of order \(n\) has the eigenvalue \(\frac nq\) with multiplicity at least \(n-1\), and \(MM^{T}\one_n=\frac{n\kappa}{q}\one_n\) by the first two identities. Therefore \(MM^{T}=\frac nqI+\frac{\kappa-1}{q}J\), and \(\frac{\kappa-1}{q}=\lambda_A\), \(\frac nq=\kappa-\lambda_A\).
\end{proof}

\subsection{Auxiliary matrices and decompositions of the complete graph}\label{sec:auxreview}
We illustrate on the smallest case how auxiliary matrices are used to decompose the complete graph~\cite{KSS2017}; \cref{sec:construction,sec:shift} generalize each step. For \(q=2\), the auxiliary matrices \(C_{y,0}\), \(C_{0,0}\), \(C_{1,0}\) of~\cite{KSS2017}, built from \(\F_2\), are
\[
C_1=I_2\otimes J_2,\qquad C_2=J_2\otimes I_2,\qquad C_3=I_2\otimes I_2+U\otimes U,
\]
where \(U\) is the \(2\times2\) permutation matrix of order \(2\), and \(J_4-C_u\) is the auxiliary matrix \(C_{a,1}\) of~\cite{KSS2017} with the same first index. These satisfy \(C_u^2=2C_u\) and \(C_uC_{u'}=J_4\) for \(u\ne u'\). Index four copies of \(\F_2^2\) by \(i\in\F_2^2\) and label the nonzero elements \((0,1),(1,0),(1,1)\) of \(\F_2^2\) by \(1,2,3\); then \(L(i,j)=i+j\) is a symmetric Latin square of order \(4\) with constant diagonal. Let \(\Gamma_0\) and \(\Gamma_1\) be the \(16\times16\) matrices with zero diagonal blocks and \((i,j)\) blocks \(C_{L(i,j)}\) and \(J_4-C_{L(i,j)}\), respectively, for \(i\ne j\), and put \(B=I_4\otimes J_4\).
\begin{enumerate}[label=(\roman*)]
\item\emph{Graphs.} The identities above give \(\Gamma_0^2=\Gamma_1^2=6I+2(J-I)\), so \(\Gamma_0\) and \(\Gamma_1\) are strongly regular graphs with parameters \((16,6,2,2)\).
\item\emph{Decomposition.} \(\Gamma_0+\Gamma_1=(J_4-I_4)\otimes J_4\), so \(K_{16}\) is the edge-disjoint union of \(\Gamma_0\), \(\Gamma_1\) and \(4K_4\).
\item\emph{Algebra.} \(\Gamma_0\Gamma_1=\Gamma_1\Gamma_0=4(B-I)+2(J-B)\). Hence \(M=\Gamma_0-\Gamma_1\) satisfies \(M^2=16(I-\frac14B)\), the matrix \(M/4\) has eigenvalues \(\pm1\), each with multiplicity \(6\), and \(0\) on the range of \(B\), and the algebra generated by \(I\) and \(M\) is \(\C\oplus\C\oplus\C\). Together with \(B-I\), the graphs \(\Gamma_0,\Gamma_1\) form an amorphic association scheme~\cite{vanDam2003}.
\end{enumerate}

\subsection{Kabanov's construction}\label{sec:kabanovdef}
We recall~\cite[Construction~3]{Kabanov2026JCD} in the case used in this paper. Let \(N\) be a symmetric incidence matrix with zero diagonal of a symmetric \((v,\kappa,\lambda)\) design. For each \(i\in[v]\) let \(\mathcal D_i\) be an affine design on a point set \(\mathcal P_i\) with \(\kappa\) parallel classes \(\mathcal B_i^1,\dots,\mathcal B_i^\kappa\) of \(q\) blocks each, all with the same parameters. Let \(L=(e(i,j))\) be a table with \(e(i,j)\in[\kappa]\) defined when \(N_{ij}=1\), such that for every \(i\) the map \(j\mapsto e(i,j)\) is a bijection from \(\{j:N_{ij}=1\}\) onto \([\kappa]\). For every \(i,j\) with \(N_{ij}=1\) let \(\varphi_{ij}\colon\mathcal B_i^{e(i,j)}\to\mathcal B_j^{e(j,i)}\) be a bijection with \(\varphi_{ji}=\varphi_{ij}^{-1}\). For \(x\in\mathcal P_i\) let \(\mathcal B_i^{u}(x)\) be the block of the \(u\)-th class containing \(x\). Kabanov's graph has vertex set \(\bigsqcup_i\mathcal P_i\), and \(x\in\mathcal P_i\), \(y\in\mathcal P_j\) are adjacent if and only if \(N_{ij}=1\) and
\[
y\in\varphi_{ij}\bigl(\mathcal B_i^{e(i,j)}(x)\bigr).
\]
Kabanov showed by counting common neighbours that this graph is a divisible design graph~\cite[Theorem~3]{Kabanov2026JCD}; he also allows \(N\) with nonzero diagonal entries when \(q\) is a power of \(2\) and the bijections \(\varphi_{ii}\) have no fixed blocks. In~\cite[Construction~1]{Kabanov2026DM} the adjacency inside the blocks with \(N_{ij}=1\) is the negation of the rule above.

\section{Auxiliary matrices}\label{sec:aux}

In this section we recall the auxiliary matrices used for the decompositions of the complete graph in~\cite{KSS2017,KharaghaniSuda2018,CIKSS2026} and identify them with parallel classes of \(\AG(2,p)\) and bijections between them. This identification is what allows the decompositions to be carried over to Kabanov's graphs in \cref{sec:construction,sec:shift}.

Let \(H=(h_{at})_{a,t\in\F_p}\) be a \(\GH(p,1)\) with rows and columns indexed by \(\F_p\), and let
\[
R_a=\bigl[U^{h_{a0}}\ U^{h_{a1}}\ \cdots\ U^{h_{a,p-1}}\bigr]\quad(a\in\F_p),\qquad R_\infty=I_p\otimes\one_p^{T}
\]
be the expansions of its rows, which are \(p\times p^2\) matrices, together with one extra \(p\times p^2\) matrix \(R_\infty\). The products \(C_a=R_a^{T}R_a\) (\(a\in\F_p\)) and \(C_\infty=R_\infty^{T}R_\infty=I_p\otimes J_p\) are the auxiliary matrices used in~\cite{CIKSS2026}. We index the columns of these matrices by pairs \((x,y)\in\F_p^2\), where \(x\) is the block and \(y\) the column inside the block, and we regard \((x,y)\) as a point of \(\AG(2,p)\).

\begin{lemma}\label{lem:gh}
The matrices \(R_a\) (\(a\in\F_p\cup\{\infty\}\)) form a \((p^2,p)\)-affine system. If \(H=H_p\), then \(R_a[b,(x,y)]=1\) if and only if \(y=ax+b\) (\(a\in\F_p\)) and \(R_\infty[c,(x,y)]=1\) if and only if \(x=c\); that is, \(R_a\) is the incidence matrix between the lines of slope \(a\) of \(\AG(2,p)\) and its points, and \(R_\infty\) is that of the vertical lines.
\end{lemma}

\begin{proof}
Each block of \(R_a\) is a permutation matrix, which gives \ref{A1}, and \ref{A2}, \ref{A3} for \(a,a'\in\F_p\) are \eqref{eq:gh} with \(\lambda=1\). For \(a\in\F_p\), \(R_\infty R_a^{T}=\sum_t e_t\one^{T}(U^{h_{at}})^{T}=\sum_t e_t\one^{T}=J_p\), where \(e_t\) is the \(t\)-th standard unit vector of length \(p\), and \(R_\infty R_\infty^{T}=pI_p\). The description for \(H=H_p\) follows from \(U^{ax}[b,y]=1\iff y=ax+b\).
\end{proof}

\begin{example}\label{ex:p3}
Let \(p=3\) and \(H=H_3\), so that \(H_3\) has rows \((0,0,0)\), \((0,1,2)\), \((0,2,1)\). Then \(R_\infty=I_3\otimes\one_3^{T}\) and the row \((0,1,2)\) gives
\[
R_\infty=\left(\begin{smallmatrix}1&1&1&0&0&0&0&0&0\\0&0&0&1&1&1&0&0&0\\0&0&0&0&0&0&1&1&1\end{smallmatrix}\right),\qquad
R_1=[I\ U\ U^2]=\left(\begin{smallmatrix}1&0&0&0&1&0&0&0&1\\0&1&0&0&0&1&1&0&0\\0&0&1&1&0&0&0&1&0\end{smallmatrix}\right).
\]
The columns are the points \((x,y)\) in the order \((0,0),(0,1),\dots,(2,2)\); the rows of \(R_\infty\) are the vertical lines \(x=c\), and the row \(b\) of \(R_1\) is the line \(y=x+b\).
\end{example}

For a point \(\xi\in\F_p^2\) and \(a\in\F_p\cup\{\infty\}\) let \(\ell_a(\xi)\) be the index of the row of \(R_a\) containing \(\xi\), that is, the line of the \(a\)-th parallel class through \(\xi\). If \(H=H_p\), then \(\ell_a(x,y)=y-ax\) for \(a\in\F_p\) and \(\ell_\infty(x,y)=x\) by \cref{lem:gh}. For a permutation \(\sigma\) of \(\F_p\) put
\[
X_{a,\sigma}=R_a^{T}P_\sigma R_a .
\]

\begin{lemma}\label{lem:X}
Let \(a,b\in\F_p\cup\{\infty\}\) and let \(\sigma,\tau\) be permutations of \(\F_p\).
\begin{enumerate}[label=(\alph*)]
\item\label{X:entry} \(X_{a,\sigma}[\xi,\xi']=1\) if and only if \(\ell_a(\xi')=\sigma(\ell_a(\xi))\); that is, \(\xi'\) lies on the image under \(\sigma\) of the line of the \(a\)-th class through \(\xi\).
\item\label{X:prod} The following identities hold:
\begin{equation}\label{eq:Xprod}
\begin{gathered}
X_{a,\sigma}^{T}=X_{a,\sigma^{-1}},\qquad X_{a,\sigma}X_{a,\tau}=p\,X_{a,\tau\circ\sigma},\\
X_{a,\sigma}X_{b,\tau}=J_{p^2}\ (a\ne b),\qquad \sum_{h\in\F_p}X_{a,\sigma\circ\tau_h}=J_{p^2},
\end{gathered}
\end{equation}
where \(\tau_h(z)=z+h\).
\end{enumerate}
\end{lemma}

\begin{proof}
\ref{X:entry} We have \(X_{a,\sigma}[\xi,\xi']=\sum_{c,d}R_a[c,\xi]\,P_\sigma[c,d]\,R_a[d,\xi']\). By \cref{lem:gh} and \ref{A1} the only nonzero term is \(c=\ell_a(\xi)\), \(d=\ell_a(\xi')\), so \(X_{a,\sigma}[\xi,\xi']=P_\sigma[\ell_a(\xi),\ell_a(\xi')]\), which is \(1\) exactly when \(\ell_a(\xi')=\sigma(\ell_a(\xi))\).

\ref{X:prod} The first identity is \(P_\sigma^{T}=P_{\sigma^{-1}}\). By \cref{lem:gh} and \ref{A2}, \(R_aR_a^{T}=pI_p\), so \(X_{a,\sigma}X_{a,\tau}=p\,R_a^{T}P_\sigma P_\tau R_a=p\,X_{a,\tau\circ\sigma}\). For \(a\ne b\), \ref{A3} gives \(R_aR_b^{T}=J_p\), so \(X_{a,\sigma}X_{b,\tau}=R_a^{T}P_\sigma JP_\tau R_b=R_a^{T}JR_b=J_{p^2}\) by \cref{lem:system}. Finally \(\sum_hP_{\sigma\circ\tau_h}=\sum_hP_{\tau_h}P_\sigma=JP_\sigma=J\), and \(R_a^{T}JR_a=J_{p^2}\).
\end{proof}

The auxiliary matrices of~\cite[Appendix~A]{CIKSS2026} for \(q=p\) are the following special cases; there \(\infty\) is denoted by \(y\). For \(a,\alpha\in\F_p\) let
\[
D_{a,\alpha}=\sum_{c\in\F_p}U^c\otimes U^{ac+\alpha},\qquad C_{a,\alpha}=\sum_{c\in\F_p}U^cR\otimes U^{ac+\alpha}R,
\]
and \(D_{\infty,\alpha}=U^\alpha\otimes J_p\), \(C_{\infty,\alpha}=U^\alpha R\otimes J_p\).

\begin{proposition}\label{prop:cikss}
Let \(H=H_p\). For \(a,\alpha\in\F_p\),
\[
D_{a,\alpha}=R_a^{T}U^{\alpha}R_a,\qquad C_{a,\alpha}=R_a^{T}U^{\alpha-a}R\,R_a,
\]
and \(D_{\infty,\alpha}=R_\infty^{T}U^\alpha R_\infty\), \(C_{\infty,\alpha}=R_\infty^{T}U^\alpha R\,R_\infty\). In particular \(D_{a,\alpha}=X_{a,\tau}\) with \(\tau(z)=z+\alpha\), and \(C_{a,\alpha}=X_{a,\psi}\) with \(\psi(z)=a-1-\alpha-z\) (\(\psi(z)=-1-\alpha-z\) for \(a=\infty\)).
\end{proposition}

\begin{proof}
Let \(\xi=(x,y)\) and \(\xi'=(x',y')\). The term \(c\) of \(D_{a,\alpha}\) has a \(1\) in position \((\xi,\xi')\) exactly when \(x'=x+c\) and \(y'=y+ac+\alpha\); eliminating \(c=x'-x\) gives \(y'-ax'=y-ax+\alpha\), i.e.\ \(\ell_a(\xi')=\ell_a(\xi)+\alpha\), and at most one \(c\) qualifies. The term \(c\) of \(C_{a,\alpha}\) has a \(1\) exactly when \(x+x'=-1-c\) and \(y+y'=-1-ac-\alpha\); eliminating \(c\) gives \(\ell_a(\xi')=-\ell_a(\xi)+a-1-\alpha\). By \cref{lem:X}\ref{X:entry}, these are the matrices \(R_a^{T}P_\sigma R_a\) for the stated \(\sigma\), and \(P_\sigma=U^\alpha\), resp.\ \(P_\sigma=U^{\alpha-a}R\). The case \(a=\infty\) is the same computation with \(\ell_\infty(x,y)=x\).
\end{proof}

\Cref{prop:cikss} explains the identities of~\cite[Lemma~2]{CIKSS2026}: \(\psi\) is an involution, so \(C_{a,\alpha}^{T}=C_{a,\alpha}\); \(\tau^{-1}(z)=z-\alpha\), so \(D_{a,\alpha}^{T}=D_{a,-\alpha}\); the products are \eqref{eq:Xprod}; and \(\sum_\alpha C_{a,\alpha}=\sum_\alpha D_{a,\alpha}=J_{p^2}\). In \cref{ex:p3}, for the translation \(\tau(z)=z+1\), the matrix \(X_{1,\tau}=R_1^{T}UR_1=D_{1,1}\) joins every point on the line \(y=x+b\) to the three points on the line \(y=x+b+1\).

\section{Adjacency matrices of Kabanov's graphs}\label{sec:construction}

The adjacency matrices of the following construction generalize the auxiliary matrices \(X_{a,\sigma}=R_a^{T}P_\sigma R_a\) of \cref{sec:aux}: the two copies of \(R_a\) are replaced by the incidence matrices of parallel classes of two possibly different affine designs, and the blocks are arranged according to a symmetric design. We show that the resulting graphs are exactly Kabanov's graphs (\cref{prop:kabanov}). This allows us to treat their products by linear algebra. As a consequence, we obtain the eigenvalues with their multiplicities (\cref{cor:spectrum}).

\begin{construction}\label{con:main}
Let the following data be given.
\begin{enumerate}[label=(D\arabic*)]
\item\label{D1} A symmetric \((0,1)\)-matrix \(N\) of order \(v\) with zero diagonal such that \(N^2=(\kappa-\lambda)I_v+\lambda J_v\).
\item\label{D2} For each \(i\in[v]\), an \((n,q)\)-affine system \(R^{(i)}_1,\dots,R^{(i)}_\kappa\), all with the same \(n\), \(q\) and \(\kappa=(n-1)/(q-1)\).
\item\label{D3} A table \(e(i,j)\in[\kappa]\), defined when \(N_{ij}=1\), such that for every \(i\) the map \(j\mapsto e(i,j)\) is a bijection from \(\{j:N_{ij}=1\}\) onto \([\kappa]\).
\item\label{D4} Permutations \(\sigma_{ij}\) of the row index set \(Q\) of the systems, defined when \(N_{ij}=1\), with \(\sigma_{ji}=\sigma_{ij}^{-1}\).
\end{enumerate}
Put \(X_{ij}=R^{(i)T}_{e(i,j)}\,P_{\sigma_{ij}}\,R^{(j)}_{e(j,i)}\) and let \(\Gamma=\Gamma(N,R,e,\sigma)\) be the \(vn\times vn\) block matrix whose \((i,j)\) block is \(N_{ij}X_{ij}\).
\end{construction}

Since \(X_{ji}=X_{ij}^{T}\) by \ref{D4} and \(N\) has zero diagonal, \(\Gamma\) is a symmetric \((0,1)\)-matrix with zero diagonal, that is, the adjacency matrix of a graph, which we also denote by \(\Gamma\). Let \(B=I_v\otimes J_n\) be the adjacency matrix of the partition into \(v\) classes of size \(n\), with loops.

\begin{theorem}\label{thm:ddg}
With \(\lambda_A=\frac{n-q}{q(q-1)}\),
\[
\Gamma^2=\frac{\kappa n}{q}I+\frac{n\lambda_A}{q}(B-I)+\frac{n\lambda}{q^2}(J-B).
\]
Hence \(\Gamma\) is a divisible design graph with parameters \(\bigl(vn,\frac{\kappa n}{q},\frac{n\lambda_A}{q},\frac{n\lambda}{q^2},v,n\bigr)\). If \(\lambda=q\lambda_A\), that is, if \(N=J_v-I_v\) with \(v=\kappa+1\), then \(\Gamma\) is a strongly regular graph with parameters \(\bigl(vn,\frac{\kappa n}{q},\frac{n\lambda}{q^2},\frac{n\lambda}{q^2}\bigr)\). Moreover \(\Gamma'=N\otimes J_n-\Gamma\) is a divisible design graph with parameters
\[
\Bigl(vn,\ \frac{\kappa n(q-1)}{q},\ \frac{\kappa n(q-2)}{q}+\frac{n\lambda_A}{q},\ \frac{n\lambda(q-1)^2}{q^2},\ v,\ n\Bigr).
\]
\end{theorem}

\begin{proof}
The \((i,k)\) block of \(\Gamma^2\) is \(\sum_jN_{ij}N_{jk}X_{ij}X_{jk}\), and
\[
X_{ij}X_{jk}=R^{(i)T}_{e(i,j)}P_{\sigma_{ij}}\bigl(R^{(j)}_{e(j,i)}R^{(j)T}_{e(j,k)}\bigr)P_{\sigma_{jk}}R^{(k)}_{e(k,j)} .
\]
If \(i\ne k\), then \(e(j,i)\ne e(j,k)\) by \ref{D3}, so the middle factor is \(\frac{n}{q^2}J_q\) by \ref{A3}, and \(R^{T}P J P R'=R^{T}JR'=J_n\) by \cref{lem:system}; the number of \(j\) with \(N_{ij}N_{jk}=1\) is \(\lambda\). If \(i=k\), the middle factor is \(\frac nqI_q\) by \ref{A2}, and \(P_{\sigma_{ij}}P_{\sigma_{ji}}=I\), so the diagonal block is \(\frac nq\sum_jN_{ij}R^{(i)T}_{e(i,j)}R^{(i)}_{e(i,j)}=\frac nq\sum_{u=1}^\kappa R^{(i)T}_uR^{(i)}_u\) by \ref{D3}, which is \(\frac nq\bigl((\kappa-\lambda_A)I_n+\lambda_AJ_n\bigr)=\frac{\kappa n}{q}I_n+\frac{n\lambda_A}{q}(J_n-I_n)\) by \cref{lem:system}. If \(\lambda=q\lambda_A\), the coefficients of \(B-I\) and \(J-B\) agree, and a graph with no edges inside the classes is then strongly regular with \(\lambda=\mu\). Since \(q\lambda_A=\kappa-1\) and \(\lambda(v-1)=\kappa(\kappa-1)\), the condition \(\lambda=q\lambda_A\) means \(v=\kappa+1\), that is, \(N=J-I\).

For \(\Gamma'\), write \(X'_{ij}=J_n-X_{ij}\); since \(X_{ij}J=JX_{ij}=\frac nqJ\), we get \(X'_{ij}X'_{jk}=n(1-\frac2q)J+X_{ij}X_{jk}\), and summing as above gives the parameters.
\end{proof}

The graphs of \cref{con:main} are exactly Kabanov's graphs of \cref{sec:kabanovdef}.

\begin{proposition}\label{prop:kabanov}
Let \(\mathcal D_i\) be the affine design whose \(u\)-th parallel class has incidence matrix \(R^{(i)}_u\), with the blocks of each class labelled by \(Q\). Then \(\Gamma(N,R,e,\sigma)\) is the graph of~\cite[Construction~3]{Kabanov2026JCD} for the collection \(\{\mathcal D_i\}\), the table \(L=(e(i,j))\) and the bijections \(\sigma_{ij}\), and \(\Gamma'=N\otimes J-\Gamma\) is the graph of~\cite[Construction~1]{Kabanov2026DM}. Conversely, every graph of these constructions with \(N\) of zero diagonal is of the form \(\Gamma(N,R,e,\sigma)\), resp.\ \(N\otimes J-\Gamma(N,R,e,\sigma)\).
\end{proposition}

\begin{proof}
For \(x\) in the \(i\)-th copy and \(y\) in the \(j\)-th copy,
\[
X_{ij}[x,y]=\sum_{c,d}R^{(i)}_{e(i,j)}[c,x]\,P_{\sigma_{ij}}[c,d]\,R^{(j)}_{e(j,i)}[d,y]=1
\]
if and only if the block of the \(e(j,i)\)-th class of \(\mathcal D_j\) containing \(y\) is the image under \(\sigma_{ij}\) of the block of the \(e(i,j)\)-th class of \(\mathcal D_i\) containing \(x\). This is Kabanov's adjacency rule; \ref{D3} is his condition on the table, and \ref{D4} is \(\sigma_{ji}=\sigma_{ij}^{-1}\). In Construction~1 of~\cite{Kabanov2026DM} the rule is replaced by its negation inside the blocks with \(N_{ij}=1\), which gives \(N\otimes J-\Gamma\). The converse holds because, by~\cite[Theorem~5.21]{Stinson2004}, the incidence matrices of the parallel classes of every affine design satisfy \ref{A1}--\ref{A3}, that is, form an affine system.
\end{proof}

\begin{corollary}\label{cor:aux}
Let \(N=J-I\) of order \(v=\kappa+1\) with \(\kappa=(p^m-1)/(p-1)\), let every \(R^{(i)}\) be the system of hyperplanes of \(\AG(m,p)\)~\cite[Section~5.2.3]{Stinson2004}, and let \(e\) be any table satisfying \ref{D3}, for instance \(e(i,j)\equiv j-i\pmod{\kappa+1}\), or, when \(\kappa+1\) is even, a symmetric Latin square of order \(\kappa+1\) with constant diagonal (the diagonal symbol deleted). Then for every choice of the \(\sigma_{ij}\) the graph \(\Gamma\) is a strongly regular graph with parameters
\[
\Bigl(p^m\frac{p^m+p-2}{p-1},\ p^{m-1}\frac{p^m-1}{p-1},\ p^{m-1}\frac{p^{m-1}-1}{p-1},\ p^{m-1}\frac{p^{m-1}-1}{p-1}\Bigr).
\]
For \(p=3\), \(m=2,3\) this gives \(\SRG(45,12,3,3)\) and \(\SRG(378,117,36,36)\).
\end{corollary}

\begin{proof}
Apply \cref{thm:ddg} with \(n=p^m\), \(q=p\) and \(\lambda=\kappa-1\). A symmetric Latin square with constant diagonal exists for every even order; since \(\kappa=1+p+\dots+p^{m-1}\), the order \(\kappa+1\) is even exactly when \(p=2\) or \(m\) is odd.
\end{proof}

The strongly regular graphs of \cref{cor:aux} have the parameters of Wallis~\cite{Wallis1971} and Fon-Der-Flaass~\cite{FDF2002}, and the divisible design graphs \(\Gamma'\) with \(N=J-I\) and \(q=3\) appear in~\cite[Table~2]{Kabanov2026JCD}; the smallest case is \cref{ex:45}. The bijections realised by the auxiliary matrices \(C_{a,\alpha},D_{a,\alpha}\) of \cref{prop:cikss} are the maps \(z\mapsto\pm z+c\); for \(p=3\) these are all the permutations of \(\F_3\). When \(e\) is a symmetric Latin square and all \(\sigma_{ij}\) are of the form \(z\mapsto c-z\), the matrix \(\Gamma\) is the incidence matrix of a symmetric design of Ionin and Kharaghani~\cite[Theorem~4.1]{IK2003}, built from a symmetric array of hyperplanes of \(\AG(m,p)\) and their translates with empty diagonal.

\begin{example}\label{ex:45}
Let \(p=3\), \(m=2\), so that every \(R^{(i)}\) is the system \(R_0,R_1,R_2,R_\infty\) of \cref{ex:p3}, let \(N=J_5-I_5\) and \(e(i,j)\equiv j-i\pmod 5\), identifying \(\{1,2,3,4\}\) with the classes \(R_0,R_1,R_2,R_\infty\). For every choice of permutations \(\sigma_{ij}\) of \(\F_3\) with \(\sigma_{ji}=\sigma_{ij}^{-1}\), the matrix \(\Gamma\) is the adjacency matrix of an \(\SRG(45,12,3,3)\), and \(\Gamma'=(J_5-I_5)\otimes J_9-\Gamma\) is a \(\DDG(45,24,15,12,5,9)\). For instance, all \(\sigma_{ij}\) may be the identity, so that every off-diagonal block of \(\Gamma\) is \(R^{T}_{e(i,j)}R_{e(j,i)}\). The construction is carried out step by step in \cref{app:srg45}.
\end{example}

\section{Decompositions of the complete graph and difference sets}\label{sec:shift}

In this section we consider several graphs of \cref{con:main} on the same vertex set, with the same \(N\), affine systems and table \(e\), and with bijections that differ by translations. For auxiliary matrices from \(\F_{2^m}\), such a family decomposes the complete graph~\cite{KSS2017}, and for \(\F_q\) with \(q\) odd the products of the graphs of such a family depend only on the differences of the shifts~\cite[Section~5]{KharaghaniSuda2018}. We show that the decomposition and the product rule hold for every graph of \cref{con:main} (\cref{thm:decomposition}), determine the spectra of the unions of the graphs of such a family, and characterize the unions that are divisible design graphs.

From now on fix the data \ref{D1}--\ref{D4} of \cref{con:main}, a total order on \([v]\), and an abelian group \(G\) of order \(q\); \(G\) need not be the additive group of \(\F_q\). Identify the row index set \(Q\) of all affine systems with \(G\) by an arbitrary bijection, and for \(g\in G\) let \(L_g\) be the translation \(c\mapsto c+g\) of \(G\). For \(g\in G\) define new bijections by
\[
\sigma^g_{ij}=L_g\circ\sigma_{ij}\quad(i<j),\qquad \sigma^g_{ji}=(\sigma^g_{ij})^{-1},
\]
that is, \(\sigma^g_{ij}\) is \(\sigma_{ij}\) followed by the translation by \(g\),
and put \(\Gamma_g=\Gamma(N,R,e,\sigma^g)\) and \(X_{ij}(g)=R^{(i)T}_{e(i,j)}P_{\sigma^g_{ij}}R^{(j)}_{e(j,i)}\). For the auxiliary matrices of order \(p^2\) and \(\sigma_{ij}\) a translation or a reflection, this replaces \(D_{a,\alpha}\) by \(D_{a,\alpha+g}\) and \(C_{a,\alpha}\) by \(C_{a,\alpha-g}\) (\cref{prop:cikss}), which is the shift used in~\cite[Appendix~A]{CIKSS2026} and~\cite[Section~5]{KharaghaniSuda2018}.

\begin{theorem}\label{thm:decomposition}
\begin{enumerate}[label=(\alph*)]
\item Each \(\Gamma_g\) is a divisible design graph with the parameters of \cref{thm:ddg}.
\item \(\sum_{g\in G}\Gamma_g=N\otimes J_n\). In particular, if \(N=J-I\), then \(K_{vn}\) is the edge-disjoint union of the \(q\) graphs \(\Gamma_g\) and \(v\) copies of \(K_n\).
\item There is a function \(F\colon G\to M_{vn}(\C)\) such that \(\Gamma_g\Gamma_h=F(h-g)\) for all \(g,h\in G\). Explicitly,
\[
F(d)=\frac{n\lambda}{q^2}(J_v-I_v)\otimes J_n+\frac nq\Delta(d),
\]
where \(\Delta(d)\) is block diagonal with \(i\)-th block
\[
\Delta_i(d)=\sum_{j<i}N_{ij}\,R^{(i)T}_{e(i,j)}P_{L_d}R^{(i)}_{e(i,j)}+\sum_{j>i}N_{ij}\,R^{(i)T}_{e(i,j)}P_{\sigma_{ij}^{-1}\circ L_{-d}\circ\sigma_{ij}}R^{(i)}_{e(i,j)} .
\]
\end{enumerate}
\end{theorem}

\begin{proof}
(a) is \cref{thm:ddg}, since the \(\sigma^g_{ij}\) satisfy \ref{D4}. (b) For fixed \(i<j\), \(\sum_gP_{L_g\circ\sigma_{ij}}=\sum_gP_{\sigma_{ij}}P_{L_g}=P_{\sigma_{ij}}J=J\), so \(\sum_gX_{ij}(g)=R^{T}JR'=J_n\) by \cref{lem:system}; the same holds for \(i>j\) by transposition.

(c) The off-diagonal blocks of \(\Gamma_g\Gamma_h\) are \(\frac{n\lambda}{q^2}J_n\) by the proof of \cref{thm:ddg}, for any bijections. The \(i\)-th diagonal block is \(\frac nq\sum_jN_{ij}R^{(i)T}_{e(i,j)}P_{\sigma^g_{ij}}P_{\sigma^h_{ji}}R^{(i)}_{e(i,j)}\) and \(P_{\sigma^g_{ij}}P_{\sigma^h_{ji}}=P_{\sigma^h_{ji}\circ\sigma^g_{ij}}\). If \(j<i\), then \(\sigma^g_{ij}=\sigma_{ji}^{-1}\circ L_{-g}\) and \(\sigma^h_{ji}=L_h\circ\sigma_{ji}\), so \(\sigma^h_{ji}\circ\sigma^g_{ij}=L_{h-g}\). If \(j>i\), then \(\sigma^g_{ij}=L_g\circ\sigma_{ij}\) and \(\sigma^h_{ji}=\sigma_{ij}^{-1}\circ L_{-h}\), so \(\sigma^h_{ji}\circ\sigma^g_{ij}=\sigma_{ij}^{-1}\circ L_{g-h}\circ\sigma_{ij}\). Both depend only on \(d=h-g\).
\end{proof}

Let \(W=\{w\in\C^{vn}:Bw=0\}\), where \(B=I_v\otimes J_n\), be the space of vectors whose entries sum to \(0\) on every class; \(\dim W=v(n-1)=(q-1)v\kappa\). Every block \(X_{ij}(g)\) has all row and column sums \(\frac nq\), so \(\Gamma_gB=B\Gamma_g\), and the symmetric matrix \(\Gamma_g\) maps \(W\) to itself. For \(g,d\in G\) put
\[
A_g=\tfrac qn\Gamma_g|_W,\qquad T_d=A_dA_0=\tfrac{q^2}{n^2}\Gamma_d\Gamma_0|_W .
\]
For a character \(\chi\) of \(G\) let \(W_\chi=\{w\in W:T_dw=\chi(d)w\text{ for all }d\in G\}\).

\begin{theorem}\label{thm:eigenspaces}
\begin{enumerate}[label=(\alph*)]
\item\label{eig:rel} Each \(A_g\) is a symmetric involution, and \(A_gA_h=T_{g-h}\) for all \(g,h\in G\). Moreover \(T_dT_e=T_{d+e}\) and \(T_dA_g=A_gT_{-d}=A_{d+g}\) for all \(d,e,g\in G\).
\item\label{eig:dim} \(W\) is the orthogonal direct sum of the spaces \(W_\chi\), \(\chi\ne1\); each of them has dimension \(v\kappa\), and \(W_1=0\). Each \(A_g\) maps \(W_\chi\) onto \(W_{\bar\chi}\).
\item\label{eig:real} For every real character \(\chi\ne1\), the restriction of \(A_0\) to \(W_\chi\) has the eigenvalues \(1\) and \(-1\), each with multiplicity \(\frac{v\kappa}{2}\).
\end{enumerate}
\end{theorem}

\begin{proof}
\ref{eig:rel} The matrices \(J-B\) and \(B\) vanish on \(W\), so \cref{thm:decomposition}(c) gives \(A_gA_h=\frac qn\Delta(h-g)|_W\), which depends only on \(g-h\). Hence \(A_gA_h=A_{g-h}A_0=T_{g-h}\). By \cref{lem:system} and \ref{D3}, \(\Delta_i(0)=\sum_{u=1}^{\kappa}R^{(i)T}_uR^{(i)}_u=\frac nqI+\lambda_AJ\), since \(\kappa-\lambda_A=\frac nq\); hence \(A_g^2=T_0=1_W\). Using \(A_0A_e=T_{-e}=A_{-e}A_0\) we get \(T_dT_e=A_dA_{-e}A_0^2=T_{d+e}\), \(T_dA_g=A_dA_{-g}A_0=T_{d+g}A_0=A_{d+g}\) and \(A_gT_{-d}=A_gA_{-d}A_0=A_{g+d}\).

\ref{eig:dim} Each \(T_d\) is a product of two symmetric involutions, hence orthogonal, and the \(T_d\) commute by~\ref{eig:rel}. So \(W\) is the orthogonal direct sum of the \(W_\chi\), and \(P_\chi=\frac1q\sum_{d}\overline{\chi(d)}T_d\) is the projection onto \(W_\chi\). We compute \(\dim W_\chi=\tr P_\chi\). We have \(T_d=\frac qn\Delta(-d)|_W\). For \(d\ne0\) every permutation in the definition of \(\Delta(-d)\) is fixed-point free, and \(\tr(R^{T}P_\pi R)=\tr(P_\pi RR^{T})=\frac nq\tr P_\pi=0\) by \ref{A2}, so \(\tr\Delta(-d)=0\). On the range of \(B\), \(\Delta(-d)\) acts as \(\frac nq\kappa\), because \(R^{T}P_\pi R\one_n=\frac nq\one_n\). Hence \(\tr T_d=\frac qn\bigl(0-v\cdot\frac nq\kappa\bigr)=-v\kappa\) for \(d\ne0\), while \(\tr T_0=\dim W=(q-1)v\kappa\). Since \(\sum_d\chi(d)\) is \(0\) for \(\chi\ne1\) and \(q\) for \(\chi=1\), we get \(\dim W_\chi=v\kappa\) for \(\chi\ne1\) and \(\dim W_1=0\). Finally, for \(w\in W_\chi\), \(T_dA_gw=A_gT_{-d}w=\overline{\chi(d)}A_gw\), so \(A_gW_\chi\subseteq W_{\bar\chi}\), with equality since \(A_g^2=1_W\).

\ref{eig:real} For real \(\chi\) the involution \(A_0\) preserves \(W_\chi\), and \(\tr(A_0P_\chi)=\frac1q\sum_d\overline{\chi(d)}\tr(A_0A_dA_0)=\frac1q\sum_d\overline{\chi(d)}\tr A_d\). On the range of \(B\), spanned by the vectors \(x\otimes\one_n\), the matrix \(\Gamma_d\) acts as \(\frac nqN\otimes\frac1nJ_n\), with trace \(\frac nq\tr N=0\); since \(\tr\Gamma_d=0\), we get \(\tr A_d=0\). So \(A_0|_{W_\chi}\) has trace \(0\).
\end{proof}

\Cref{thm:eigenspaces} determines the spectrum of every graph built from the family. For \(S\subseteq G\) put \(\Gamma_S=\sum_{g\in S}\Gamma_g\) and \(\chi(S)=\sum_{g\in S}\chi(g)\).

\begin{corollary}\label{cor:spectrum}
Let \(\theta_1,\dots,\theta_v\) be the eigenvalues of \(N\). The eigenvalues of \(\Gamma_S\), with multiplicities, are \(\frac{n|S|}{q}\theta_1,\dots,\frac{n|S|}{q}\theta_v\), together with \(\pm\frac nq|\chi(S)|\), each with multiplicity \(\frac{v\kappa}2\), for each nontrivial character \(\chi\) of \(G\). In particular every graph \(\Gamma\) of \cref{con:main} has the eigenvalues \(\frac nq\theta_1,\dots,\frac nq\theta_v\) and \(\pm\frac nq\), each of the latter with multiplicity \(\frac{v(n-1)}{2}\).
\end{corollary}

\begin{proof}
On the range of \(B\), \(\Gamma_S\) acts as \(\frac{n|S|}{q}N\otimes\frac1nJ_n\), with eigenvalues \(\frac{n|S|}{q}\theta_i\). On \(W\), \(\Gamma_S=\frac nqA_S\) with \(A_S=\sum_{g\in S}A_g\), and \(A_S^2=\sum_{g,h\in S}T_{g-h}\) acts on \(W_\chi\) as \(|\chi(S)|^2\) by \cref{thm:eigenspaces}\ref{eig:rel}. Let \(\chi\ne\bar\chi\). The symmetric matrix \(A_S\) preserves \(W_\chi\oplus W_{\bar\chi}\) and interchanges the two summands, so it anticommutes with the involution that is \(1\) on \(W_\chi\) and \(-1\) on \(W_{\bar\chi}\); hence its eigenvalues there are \(\pm|\chi(S)|\), each with multiplicity \(v\kappa\), that is, \(\frac{v\kappa}2\) for each of \(\chi\) and \(\bar\chi\). Let \(\chi\ne1\) be real. On \(W_\chi\) we have \(A_gA_0=T_g=\chi(g)\), so \(A_g=\chi(g)A_0\) and \(A_S=\chi(S)A_0\), with eigenvalues \(\pm\chi(S)\), each with multiplicity \(\frac{v\kappa}2\) by \cref{thm:eigenspaces}\ref{eig:real}. For \(S=\{0\}\) we have \(|\chi(S)|=1\) and \((q-1)\frac{v\kappa}{2}=\frac{v(n-1)}{2}\).
\end{proof}

The spectrum tells which unions are again divisible design graphs. There is much work in which difference sets are used to construct designs and graphs, see for instance~\cite{KT2003,MR2026,DGMX2026,DGX2026}; in our setting the difference-set condition is also necessary. A subset \(S\subseteq G\) of size \(s\) is a \((q,s,\mu)\) \emph{difference set} if every nonzero element of \(G\) is a difference of two elements of \(S\) in exactly \(\mu\) ways; we include the trivial cases \(s\in\{1,q-1,q\}\).

\begin{proposition}\label{prop:diffset}
Let \(S\subseteq G\) be nonempty with \(s=|S|\). Then \(\Gamma_S\) is a divisible design graph if and only if \(S\) is a \((q,s,\mu)\) difference set in \(G\). In this case
\[
\Gamma_S^2=\frac{s\kappa n}{q}I+\frac nq\bigl((s-\mu)\lambda_A+\mu\kappa\bigr)(B-I)+\frac{s^2n\lambda}{q^2}(J-B),
\]
so \(\Gamma_S\) has parameters \(\bigl(vn,\frac{s\kappa n}{q},\frac nq((s-\mu)\lambda_A+\mu\kappa),\frac{s^2n\lambda}{q^2},v,n\bigr)\).
\end{proposition}

\begin{proof}
Suppose that \(S\) is a difference set. By \cref{thm:decomposition}(c),
\[
\Gamma_S^2=\sum_{g,h\in S}F(h-g)=\frac{s^2n\lambda}{q^2}(J_v-I_v)\otimes J_n+\frac nq\Bigl(s\Delta(0)+\mu\sum_{d\ne0}\Delta(d)\Bigr).
\]
Since \(\sum_{d\in G}P_{L_d}=\sum_{d\in G}P_{\sigma^{-1}\circ L_{-d}\circ\sigma}=J\) for every permutation \(\sigma\), we have \(\sum_{d\in G}\Delta_i(d)=\sum_jN_{ij}R^{(i)T}_{e(i,j)}JR^{(i)}_{e(i,j)}=\kappa J_n\), and \(\Delta_i(0)=\frac nqI+\lambda_AJ_n\) as in the proof of \cref{thm:eigenspaces}. Hence \(s\Delta(0)+\mu\sum_{d\ne0}\Delta(d)=(s-\mu)(\frac nqI+\lambda_AB)+\mu\kappa B\), which gives the stated formula, since the coefficient of \(I\) is \(\frac{n^2(s-\mu)}{q^2}=\frac{s\kappa n}{q}-\frac nq\bigl((s-\mu)\lambda_A+\mu\kappa\bigr)\) because \(\frac nq+\lambda_A=\kappa\).

Conversely, if \(\Gamma_S\) is a divisible design graph, then \(\Gamma_S^2\) is a linear combination of \(I\), \(B\) and \(J\), and therefore acts on \(W\) as a scalar. By \cref{cor:spectrum}, the eigenvalues of \(\Gamma_S^2\) on \(W\) are \(\frac{n^2}{q^2}|\chi(S)|^2\) for all nontrivial characters \(\chi\) of \(G\), each with positive multiplicity. So \(|\chi(S)|^2=c\) for all \(\chi\ne1\). By Parseval's identity \(\sum_\chi|\chi(S)|^2=qs\), we get \(s^2+(q-1)c=qs\), that is, \(c=s-\mu\) with \(\mu=\frac{s(s-1)}{q-1}\). Hence \(S\) is a \((q,s,\mu)\) difference set by~\cite[p.~319]{Turyn1965}.
\end{proof}

\begin{example}\label{ex:441}
For \(\AG(2,7)\), \(N=J_9-I_9\), \(G=\mathbb Z_7\) and the difference set \(S=\{1,2,4\}\), \cref{prop:diffset} gives a divisible design graph with parameters \((441,168,70,63,9,49)\) for every choice of the table and the bijections. It is worked out in \cref{app:441}.
\end{example}

Taking for \(\mathcal D_i\) the point-hyperplane designs of \(\AG(m,q)\), \cref{prop:diffset} gives the following families. Here \(N\) may be \(J-I\) or the adjacency matrix of any strongly regular graph with \(\lambda=\mu\) and valency \(\kappa\). The graphs are not among Kabanov's graphs, whose valencies are \(\frac{\kappa n}{q}\) and \(\frac{\kappa n(q-1)}{q}\), and we have not found their parameters among the known constructions of divisible design graphs~\cite{HKM2011,Kabanov2022,Kabanov2026DM,Kabanov2026JCD,PS2022,MR2026,DGMX2026,DGX2026}.

\begin{corollary}\label{cor:families}
Let \(q\) be a prime power, \(m\ge2\), \(\kappa=\frac{q^m-1}{q-1}\) and \(\lambda_A=\frac{q^{m-1}-1}{q-1}\). Let \(N\) be a symmetric \((0,1)\)-matrix of order \(v\) with zero diagonal and \(N^2=(\kappa-\lambda)I+\lambda J\), that is, \(N=J_{\kappa+1}-I\) or the adjacency matrix of an \(\SRG(v,\kappa,\lambda,\lambda)\), and let \(S\) be a \((q,s,\mu)\) difference set in an abelian group of order \(q\). Then there is a divisible design graph with parameters
\[
\bigl(vq^m,\ s\kappa q^{m-1},\ q^{m-1}((s-\mu)\lambda_A+\mu\kappa),\ s^2q^{m-2}\lambda,\ v,\ q^m\bigr).
\]
\end{corollary}

\begin{proof}
Apply \cref{prop:diffset} with \(n=q^m\).
\end{proof}

The following known difference sets in abelian groups of prime power order \(q\) give examples, together with their complements \((q,q-s,q-2s+\mu)\): the Paley difference sets \(\bigl(q,\frac{q-1}2,\frac{q-3}4\bigr)\) for \(q\equiv3\pmod4\)~\cite[Theorem~3.21]{Stinson2004}; the Singer difference sets \(\bigl(q,\frac{r^d-1}{r-1},\frac{r^{d-1}-1}{r-1}\bigr)\) for \(q=\frac{r^{d+1}-1}{r-1}\) with \(r\) a prime power and \(d\ge2\)~\cite[Theorem~3.28]{Stinson2004}; the Menon difference sets \(\bigl(q,\frac{q-\sqrt q}{2},\frac{q-2\sqrt q}{4}\bigr)\) in elementary abelian groups for \(q=4^t\), \(t\ge2\), obtained from the bent function \(x_1x_2+\dots+x_{2t-1}x_{2t}\)~\cite[Theorem~4.48, Corollary~4.45]{Stinson2004}; the biquadratic residues \(\bigl(q,\frac{q-1}4,\frac{q-5}{16}\bigr)\) for primes \(q=4x^2+1\), and with zero \(\bigl(q,\frac{q+3}4,\frac{q+3}{16}\bigr)\) for primes \(q=4x^2+9\), \(x\) odd; and the octic residues \(\bigl(q,\frac{q-1}8,\frac{q-9}{64}\bigr)\) for primes \(q=8a^2+1=64b^2+9\) with \(a,b\) odd, and with zero \(\bigl(q,\frac{q+7}8,\frac{q+7}{64}\bigr)\) for primes \(q=8a^2+49=64b^2+441\) with \(a\) odd and \(b\) even~\cite{Lehmer1953}, see also~\cite{Xia2017}. Besides \(N=J-I\), strongly regular graphs with \(\lambda=\mu\) and valency \(\kappa\) occur; for instance, for \(m=2\) and a Mersenne prime \(q=2^{2d-1}-1\) the symplectic graph on \(\F_2^{2d}\setminus\{0\}\), an \(\SRG(2^{2d}-1,2^{2d-1},2^{2d-2},2^{2d-2})\)~\cite{BvM2022}, has valency \(\kappa=q+1\). \Cref{tab:ddg} lists the graphs of \cref{cor:families} with at most \(10^5\) vertices obtained from these difference sets and from \(N=J-I\) or a strongly regular graph with \(\lambda=\mu\) in the table of Brouwer~\cite{BrouwerSRG}. We checked the parameters of the smallest cases \(441\), \(735\) and \(1573\) also by computer.

\input{table_ddg}

\section*{Acknowledgments}
The work of Keita Ishizuka is supported by JSPS KAKENHI Grant Numbers JP25KK0004 and JP25K17290. Hadi Kharaghani is supported by an NSERC Discovery Grant.

\appendix
\section{\(\SRG(45,12,3,3)\) by \cref{cor:aux}, step by step}\label{app:srg45}

The construction is the usual one with auxiliary matrices, with exactly two changes:
\begin{enumerate}
\item[(i)] \emph{Twist.} The block \(R_a^{T}R_a\) is replaced by \(R_a^{T}PR_b\), where \(P\) is a \(3\times3\) permutation matrix.
\item[(ii)] \emph{Free placement.} Which matrices go into the \((i,j)\) block is chosen row by row by a table \(e(i,j)\); the only rule is that every row uses each of the four classes exactly once.
\end{enumerate}
Everything else is the same.

\subsection*{Step 1: the four matrices \(R_0,R_1,R_2,R_\infty\)}
Let \(U\) be the \(3\times3\) shift matrix, \(U[c,c+1]=1\) (indices mod \(3\)), and expand the rows of the generalized Hadamard matrix
\[
H_3=\begin{pmatrix}0&0&0\\0&1&2\\0&2&1\end{pmatrix}
\]
(row \(a\) gives \(R_a\)), together with one extra matrix \(R_\infty\):
\[
R_0=[I\ I\ I],\quad R_1=[I\ U\ U^2],\quad R_2=[I\ U^2\ U],\quad R_\infty=I_3\otimes(1\ 1\ 1).
\]
Explicitly,
\[
\setlength\arraycolsep{2.5pt}
R_0=\begin{pmatrix}1&0&0&1&0&0&1&0&0\\
0&1&0&0&1&0&0&1&0\\
0&0&1&0&0&1&0&0&1\end{pmatrix},\quad
R_1=\begin{pmatrix}1&0&0&0&1&0&0&0&1\\
0&1&0&0&0&1&1&0&0\\
0&0&1&1&0&0&0&1&0\end{pmatrix},
\]
\[
\setlength\arraycolsep{2.5pt}
R_2=\begin{pmatrix}1&0&0&0&0&1&0&1&0\\
0&1&0&1&0&0&0&0&1\\
0&0&1&0&1&0&1&0&0\end{pmatrix},\quad
R_\infty=\begin{pmatrix}1&1&1&0&0&0&0&0&0\\
0&0&0&1&1&1&0&0&0\\
0&0&0&0&0&0&1&1&1\end{pmatrix}.
\]
They satisfy
\[
R_aR_a^{T}=3I_3,\qquad R_aR_b^{T}=J_3\ (a\neq b),\qquad \sum_{a}R_a^{T}R_a=3I_9+J_9 .
\]
(The products \(R_a^{T}R_a\) are the auxiliary matrices.)

\subsection*{Step 2: the table \(e(i,j)\) (change (ii))}
Take \(N=J_5-I_5\), so there are \(5\) block rows and columns and \(4\) off-diagonal blocks in each block row. Label the classes \(1,2,3,4\leftrightarrow R_0,R_1,R_2,R_\infty\) and put \(e(i,j)\equiv j-i\pmod 5\):
\[
e=\begin{pmatrix}
-&1&2&3&4\\
4&-&1&2&3\\
3&4&-&1&2\\
2&3&4&-&1\\
1&2&3&4&-
\end{pmatrix}.
\]
Each row contains \(1,2,3,4\) once. The table need not be symmetric: here \(e(i,j)\ne e(j,i)\).

\subsection*{Step 3: the twists \(P_{ij}\) (change (i))}
For \(i<j\) put \(P_{ij}=U\), and \(P_{ji}=P_{ij}^{T}=U^2\). (Any permutation matrices with \(P_{ji}=P_{ij}^{T}\) work, including \(P_{ij}=I\).)

\subsection*{Step 4: the blocks}
For \(i\ne j\) let
\[
X_{ij}=R_{e(i,j)}^{T}\,P_{ij}\,R_{e(j,i)} ,
\]
the left factor chosen by row \(i\) of the table, the right factor by row \(j\). For example,
\[
X_{12}=R_0^{T}UR_\infty,\quad X_{21}=R_\infty^{T}U^2R_0=X_{12}^{T},\]
\[X_{13}=R_1^{T}UR_2,\quad X_{15}=R_\infty^{T}UR_0 .
\]
Since \(P_{ji}=P_{ij}^{T}\), we have \(X_{ji}=X_{ij}^{T}\). Put
\[
\Gamma=\begin{pmatrix}
0&X_{12}&X_{13}&X_{14}&X_{15}\\
X_{21}&0&X_{23}&X_{24}&X_{25}\\
X_{31}&X_{32}&0&X_{34}&X_{35}\\
X_{41}&X_{42}&X_{43}&0&X_{45}\\
X_{51}&X_{52}&X_{53}&X_{54}&0
\end{pmatrix}.
\]
This is a symmetric \((0,1)\)-matrix of order \(45\) with zero diagonal and row sum \(4\cdot3=12\). It is shown in full in \cref{fig:g45}.

\subsection*{Step 5: \(\Gamma^2=12I+3(J-I)\)}
\emph{Diagonal blocks.} The \((i,i)\) block of \(\Gamma^2\) is
\[
\sum_{j\ne i}X_{ij}X_{ji}=\sum_{j\ne i} R_{e(i,j)}^{T}P_{ij}\bigl(R_{e(j,i)}R_{e(j,i)}^{T}\bigr)P_{ij}^{T}R_{e(i,j)}=3\sum_{j\ne i} R_{e(i,j)}^{T}R_{e(i,j)},
\]
because \(RR^{T}=3I\) and \(PP^{T}=I\). By (ii), \(e(i,\cdot)\) runs over all four classes, so this equals \(3(3I_9+J_9)=12I_9+3(J_9-I_9)\).

\emph{Off-diagonal blocks.} For \(i\ne k\), the \((i,k)\) block is \(\sum_{j\ne i,k}X_{ij}X_{jk}\), and the middle factor of \(X_{ij}X_{jk}\) is \(R_{e(j,i)}R_{e(j,k)}^{T}=J_3\), since \(e(j,i)\ne e(j,k)\) by (ii). As \(PJ_3P'=J_3\) for permutation matrices \(P,P'\), each term is \(R^{T}J_3R'=J_9\). There are \(3\) such \(j\), so the block is \(3J_9\).

\begin{figure}[ht]
\centering
\includegraphics[width=0.75\textwidth]{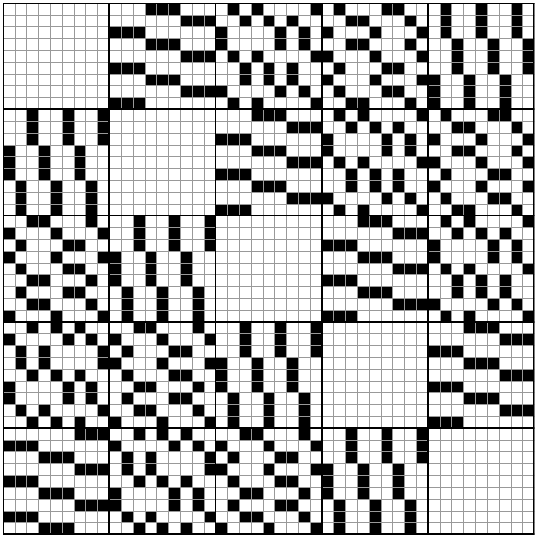}
\caption{The matrix \(\Gamma\) of Step 4 of \cref{app:srg45}. A black cell is an edge, a white cell a non-edge; the thick lines separate the five classes of size \(9\). In each class, rows and columns are the points \((x,y)\), ordered \((0,0),(0,1),\dots,(2,2)\).}
\label{fig:g45}
\end{figure}

Hence \(\Gamma^2=12I+3(J-I)\), and \(\Gamma\) is an \(\SRG(45,12,3,3)\). The twist cancels in both cases (\(PP^{T}=I\), \(PJP'=J\)), and the free placement is used only through ``each row of \(e\) uses every class once''. The matrix of \cref{fig:g45} was also checked by computer.

\subsection*{Step 6: other choices of \(e\) and \(P\)}
The proof in Step 5 uses only that every row of \(e\) contains each class once and that \(P_{ji}=P_{ij}^{T}\). Hence \emph{every} such choice gives an \(\SRG(45,12,3,3)\): there are \((4!)^5\) tables \(e\) and \(6^{10}\) choices of the \(P_{ij}\) (\(i<j\)). Changing \(e\) or \(P\) changes the matrix \(\Gamma\), and in general the graph as well. For instance, keep \(e\) of Step~2 and put \(P_{ij}=I\) for all \(i<j\), except possibly \(P_{13}\):
\begin{center}
\begin{tabular}{lc}
choice & \(|\mathrm{Aut}(\Gamma)|\)\\\hline
\(P_{13}=I\) & \(10\)\\
\(P_{13}=\left(\begin{smallmatrix}1&0&0\\0&0&1\\0&1&0\end{smallmatrix}\right)\) (\(z\mapsto -z\)) & \(18\)\\[6pt]
\(P_{13}=\left(\begin{smallmatrix}0&1&0\\1&0&0\\0&0&1\end{smallmatrix}\right)\) (\(z\mapsto 1-z\)) & \(9\)
\end{tabular}
\end{center}
So these give three pairwise non-isomorphic \(\SRG(45,12,3,3)\) (computed with nauty in SageMath). The graph of Step~4 (\(P_{ij}=U\)) is isomorphic to the first one.

\section{Shifted families and difference sets, two examples}\label{app:shift}

\Cref{sec:shift} adds one ingredient: all bijections \(\sigma_{ij}\) (\(i<j\)) are followed by the same translation \(L_g\colon z\mapsto z+g\). Since \(P_{L_g\circ\sigma}=P_\sigma U^g\), in the language of \cref{app:srg45} this replaces the twist \(P_{ij}\) by
\[
P_{ij}(g)=P_{ij}U^g\quad(i<j),\qquad P_{ji}(g)=P_{ij}(g)^{T},
\]
and \(\Gamma_g\) is the matrix of Step~4 with these twists. \Cref{app:shift45} shows the family itself on \(45\) vertices, and \cref{app:441} shows how a difference set gives a new divisible design graph.

\subsection{The family \(\Gamma_0,\Gamma_1,\Gamma_2\) on \(45\) vertices}\label{app:shift45}
\begin{figure}[ht]
\centering
\includegraphics[width=0.75\textwidth]{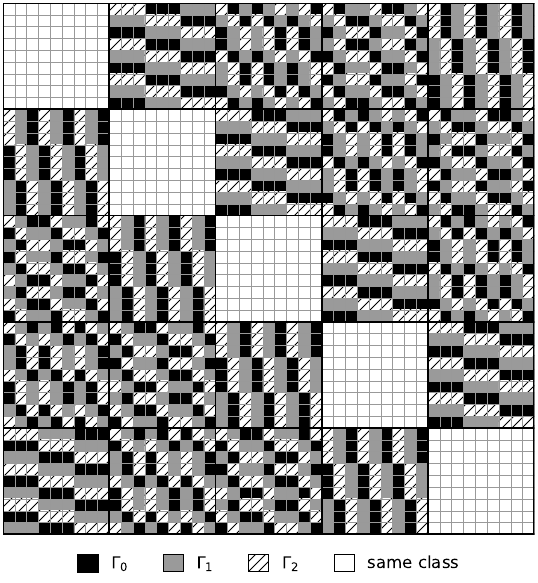}
\caption{The decomposition \(K_{45}=\Gamma_0\cup\Gamma_1\cup\Gamma_2\cup5K_9\) of \cref{app:shift45}. A black, grey or hatched cell is an edge of \(\Gamma_0\), \(\Gamma_1\) or \(\Gamma_2\), respectively; the white diagonal blocks are the pairs inside a class. Rows and columns are ordered as in \cref{fig:g45}, and the cells of \(\Gamma_0\) form the matrix of \cref{fig:g45}.}
\label{fig:labels}
\end{figure}

Take the data of \cref{app:srg45}, \(P_{ij}=U\) for \(i<j\), and \(G=\mathbb Z_3\). Then \(P_{ij}(g)=U^{1+g}\): \(\Gamma_0\) is the graph \(\Gamma\) of Step~4, \(\Gamma_1\) has \(P_{ij}=U^2\), and \(\Gamma_2\) has \(P_{ij}=I\).
\begin{enumerate}[label=(\roman*)]
\item\emph{Graphs (\cref{thm:decomposition}(a)).} Each \(\Gamma_g\) is an \(\SRG(45,12,3,3)\) by Step~5, which works for every twist.
\item\emph{Decomposition (\cref{thm:decomposition}(b)).} For \(i<j\),
\[
X_{ij}(0)+X_{ij}(1)+X_{ij}(2)=R^{T}U(I+U+U^2)R'=R^{T}J_3R'=J_9 .
\]
Hence \(\Gamma_0+\Gamma_1+\Gamma_2=(J_5-I_5)\otimes J_9\): every pair of vertices in different classes is an edge of exactly one \(\Gamma_g\), and \(K_{45}\) is the edge-disjoint union of \(\Gamma_0,\Gamma_1,\Gamma_2\) and \(5K_9\); see \cref{fig:labels}.
\item\emph{Products (\cref{thm:decomposition}(c)).} The off-diagonal blocks of \(\Gamma_g\Gamma_h\) are \(3J_9\), exactly as in Step~5. The \((i,i)\) block is, as in Step~5,
\[
\sum_{j\ne i}X_{ij}(g)X_{ji}(h)=3\sum_{j\ne i}R_{e(i,j)}^{T}\,P_{ij}(g)P_{ij}(h)^{T}\,R_{e(i,j)},
\]
but now \(P_{ij}(g)P_{ij}(h)^{T}=U^{g-h}\) for \(j>i\) and \(U^{h-g}\) for \(j<i\) instead of \(I\). It depends only on \(d=h-g\). Hence \(\Gamma_g\Gamma_h=F(h-g)\):
\begin{align*}
\Gamma_0\Gamma_1&=\Gamma_1\Gamma_2=\Gamma_2\Gamma_0=F(1),\\
\Gamma_1\Gamma_0&=\Gamma_2\Gamma_1=\Gamma_0\Gamma_2=F(2)=F(1)^{T},
\end{align*}
and \(F(1)\ne F(2)\), so \(\Gamma_0\Gamma_1\ne\Gamma_1\Gamma_0\).
\item\emph{Spectrum (\cref{cor:spectrum}).} \(N\) has the eigenvalues \(4\) and \(-1\) with multiplicity \(4\), and \(n/q=3\), so each \(\Gamma_g\) has the eigenvalues \(12\), \(-3\) with multiplicity \(4\), and \(\pm3\), each with multiplicity \(20\); in total \(12^{(1)},3^{(20)},(-3)^{(24)}\), the spectrum of an \(\SRG(45,12,3,3)\).
\item\emph{A union (\cref{prop:diffset}, trivial case).} Every \(2\)-subset of \(\mathbb Z_3\) is a \((3,2,1)\) difference set. For \(S=\{0,1\}\) the blocks of \(\Gamma_S=\Gamma_0+\Gamma_1\) are \(R^{T}MR'\) with \(M=U(I+U)\) and \(MM^{T}=I+J_3\), and the computation of \cref{app:441} gives \(\Gamma_S^2=24I+15(B-I)+12(J-B)\), a \(\DDG(45,24,15,12,5,9)\). Since \(\Gamma_S=(J_5-I_5)\otimes J_9-\Gamma_2\), this is the graph \(\Gamma'\) of \cref{thm:ddg}. A new graph needs a nontrivial difference set, and the smallest one occurs for \(q=7\).
\end{enumerate}

\subsection{The divisible design graph of \cref{ex:441}}\label{app:441}
\emph{Data.} \(p=7\), so that \(n=49\), \(q=7\), \(\kappa=8\), \(\lambda_A=1\).
\begin{itemize}
\item \(R_a=[U^{0}\ U^{a}\ U^{2a}\ \cdots\ U^{6a}]\) (\(a\in\F_7\)) and \(R_\infty=I_7\otimes(1\ \cdots\ 1)\), the eight parallel classes of \(\AG(2,7)\) (\cref{lem:gh}), labelled \(1,\dots,8\). Row \(b\) of \(R_a\) is the line \(y=ax+b\).
\item \(N=J_9-I_9\), so \(v=9\), \(\lambda=7\), and \(e(i,j)\equiv j-i\pmod 9\).
\item \(P_{ij}=P_{\sigma_{ij}}\) with \(\sigma_{ij}(z)=3z\) for \(i<j\); any permutation matrices would do.
\item \(G=\mathbb Z_7\) and \(S=\{1,2,4\}\), the nonzero squares, a \((7,3,1)\) difference set.
\end{itemize}

\emph{The blocks.} For \(i<j\) the block of \(\Gamma_S=\Gamma_1+\Gamma_2+\Gamma_4\) is
\[
\sum_{g\in S}R_{e(i,j)}^{T}P_{ij}U^gR_{e(j,i)}=R_{e(i,j)}^{T}\,M_{ij}\,R_{e(j,i)},
\]
where \(M_{ij}=P_{ij}D\) and \(D=U+U^2+U^4\) is the \(7\times7\) circulant of the difference set; \(M_{ji}=M_{ij}^{T}\). So \(\Gamma_S\) is the matrix of Step~4 of \cref{app:srg45} with the twist \(P_{ij}\) replaced by the \((0,1)\)-matrix \(M_{ij}\). In terms of points, \(\xi\) in the \(i\)-th class and \(\eta\) in the \(j\)-th class are adjacent if and only if
\[
\ell_{e(j,i)}(\eta)-\sigma_{ij}\bigl(\ell_{e(i,j)}(\xi)\bigr)\in S .
\]
Each vertex has \(8\cdot3\cdot7=168\) neighbours.

\emph{The only property of \(S\) that is used.} Since \(S\) is a \((7,3,1)\) difference set, \(DD^{T}=2I_7+J_7\), and hence \(M_{ij}M_{ij}^{T}=2I_7+J_7\) for all \(i\ne j\). For \(S=\{0,1,2\}\), which is not a difference set, the first row of \(DD^{T}\) is \((3,2,1,0,0,1,2)\), which is not of the form \(aI+bJ\).

\emph{\(\Gamma_S^2\).} This is Step~5 of \cref{app:srg45} with \(MM^{T}=2I+J\) in place of \(PP^{T}=I\). The \((i,i)\) block is
\begin{align*}
&\sum_{j\ne i}R_{e(i,j)}^{T}M_{ij}\bigl(R_{e(j,i)}R_{e(j,i)}^{T}\bigr)M_{ij}^{T}R_{e(i,j)}\\
&=7\sum_{j\ne i}R_{e(i,j)}^{T}(2I_7+J_7)R_{e(i,j)}\\
&=7\Bigl(2\sum_{u=1}^{8}R_u^{T}R_u+8J_{49}\Bigr)\\
&=7\bigl(2(7I+J)+8J\bigr)=168I+70(J-I),
\end{align*}
by \(RR^{T}=7I_7\), \(R^{T}J_7R=J_{49}\) and \cref{lem:system}. For \(i\ne k\), the middle factor of \(X^S_{ij}X^S_{jk}\) is \(R_{e(j,i)}R_{e(j,k)}^{T}=J_7\), and \(M_{ij}J_7M_{jk}=9J_7\), since every \(M_{ij}\) has row and column sums \(3\). So each of the \(\lambda=7\) terms is \(9J_{49}\), and the block is \(63J_{49}\). Hence
\[
\Gamma_S^2=168I+70(B-I)+63(J-B),
\]
and \(\Gamma_S\) is a \(\DDG(441,168,70,63,9,49)\), in accordance with \cref{prop:diffset}. As in Step~6, the table \(e\) and the \(P_{ij}\) entered only through ``each row of \(e\) uses every class once'' and the identities \(PJP'=J\) and \(P(2I+J)P^{T}=2I+J\), so every choice gives the same parameters.

\emph{Spectrum (\cref{cor:spectrum}).} \(168\), \(-21\) with multiplicity \(8\), and \(\pm7\sqrt2\), each with multiplicity \(216\), since \(|\chi(S)|^2=s-\mu=2\) for every nontrivial character \(\chi\) of \(\mathbb Z_7\).

\emph{Without a difference set.} For \(S=\{0,1,2\}\), the matrix \(\Gamma_S^2\) takes the \(15\) values \(21,28,\dots,119\) on the pairs of distinct vertices in the same class, so \(\Gamma_S\) is not a divisible design graph, as \cref{prop:diffset} predicts.

We checked all statements of this appendix by computer, for the data above and for three random choices of the table and the bijections.

\end{document}

%% file: table_ddg.tex
\begin{table}[!htbp]
\centering\footnotesize\setlength{\tabcolsep}{3.5pt}
\caption{Divisible design graphs \(\Gamma_S\) of \cref{cor:families} with at most \(10^5\) vertices. Here \(\mathcal D_i\) is the design of points and hyperplanes of \(\AG(m',q)\), \(N\) is \(J-I\) or an \(\SRG(v,\kappa,\lambda,\lambda)\) given by its parameters, and \(S\) is a \((q,s,\mu)\) difference set: P Paley, S Singer, M Menon, B biquadratic residues, B\(_0\) biquadratic residues with zero; c the complement of such a set.}\label{tab:ddg}
\begin{tabular}{rrrrrrrrlll}
\toprule
\(v\) & \(k\) & \(\lambda_1\) & \(\lambda_2\) & \(m\) & \(n\) & \(q\) & \(m'\) & \(N\) & \((q,s,\mu)\) & type\\
\midrule
441 & 168 & 70 & 63 & 9 & 49 & 7 & 2 & $K_{9}$ & $(7,3,1)$ & P,S \\
441 & 224 & 126 & 112 & 9 & 49 & 7 & 2 & $K_{9}$ & $(7,4,2)$ & P,S, c \\
735 & 168 & 70 & 36 & 15 & 49 & 7 & 2 & $(15,8,4,4)$ & $(7,3,1)$ & P,S \\
735 & 224 & 126 & 64 & 15 & 49 & 7 & 2 & $(15,8,4,4)$ & $(7,4,2)$ & P,S, c \\
1573 & 660 & 297 & 275 & 13 & 121 & 11 & 2 & $K_{13}$ & $(11,5,2)$ & P \\
1573 & 792 & 429 & 396 & 13 & 121 & 11 & 2 & $K_{13}$ & $(11,6,3)$ & P, c \\
2535 & 728 & 221 & 208 & 15 & 169 & 13 & 2 & $K_{15}$ & $(13,4,1)$ & B$_0$,S \\
2535 & 1638 & 1131 & 1053 & 15 & 169 & 13 & 2 & $K_{15}$ & $(13,9,6)$ & B$_0$,S, c \\
4608 & 1632 & 608 & 576 & 18 & 256 & 16 & 2 & $K_{18}$ & $(16,6,2)$ & M \\
4608 & 2720 & 1696 & 1600 & 18 & 256 & 16 & 2 & $K_{18}$ & $(16,10,6)$ & M, c \\
5445 & 660 & 297 & 75 & 45 & 121 & 11 & 2 & $(45,12,3,3)$ & $(11,5,2)$ & P \\
5445 & 792 & 429 & 108 & 45 & 121 & 11 & 2 & $(45,12,3,3)$ & $(11,6,3)$ & P, c \\
7581 & 3420 & 1615 & 1539 & 21 & 361 & 19 & 2 & $K_{21}$ & $(19,9,4)$ & P \\
7581 & 3800 & 1995 & 1900 & 21 & 361 & 19 & 2 & $K_{21}$ & $(19,10,5)$ & P, c \\
13225 & 6072 & 2898 & 2783 & 25 & 529 & 23 & 2 & $K_{25}$ & $(23,11,5)$ & P \\
13225 & 6624 & 3450 & 3312 & 25 & 529 & 23 & 2 & $K_{25}$ & $(23,12,6)$ & P, c \\
19894 & 8379 & 3577 & 3528 & 58 & 343 & 7 & 3 & $K_{58}$ & $(7,3,1)$ & P,S \\
19894 & 11172 & 6370 & 6272 & 58 & 343 & 7 & 3 & $K_{58}$ & $(7,4,2)$ & P,S, c \\
21141 & 9828 & 4725 & 4563 & 29 & 729 & 27 & 2 & $K_{29}$ & $(27,13,6)$ & P \\
21141 & 10584 & 5481 & 5292 & 29 & 729 & 27 & 2 & $K_{29}$ & $(27,14,7)$ & P, c \\
31713 & 5952 & 1147 & 1116 & 33 & 961 & 31 & 2 & $K_{33}$ & $(31,6,1)$ & S \\
31713 & 14880 & 7192 & 6975 & 33 & 961 & 31 & 2 & $K_{33}$ & $(31,15,7)$ & P,S \\
31713 & 15872 & 8184 & 7936 & 33 & 961 & 31 & 2 & $K_{33}$ & $(31,16,8)$ & P,S, c \\
31713 & 24800 & 19995 & 19375 & 33 & 961 & 31 & 2 & $K_{33}$ & $(31,25,20)$ & S, c \\
34656 & 3420 & 1615 & 324 & 96 & 361 & 19 & 2 & $(96,20,4,4)$ & $(19,9,4)$ & P \\
34656 & 3800 & 1995 & 400 & 96 & 361 & 19 & 2 & $(96,20,4,4)$ & $(19,10,5)$ & P, c \\
46656 & 9828 & 4725 & 2028 & 64 & 729 & 27 & 2 & $(64,28,12,12)$ & $(27,13,6)$ & P \\
46656 & 10584 & 5481 & 2352 & 64 & 729 & 27 & 2 & $(64,28,12,12)$ & $(27,14,7)$ & P, c \\
53391 & 12654 & 3071 & 2997 & 39 & 1369 & 37 & 2 & $K_{39}$ & $(37,9,2)$ & B \\
53391 & 39368 & 29785 & 29008 & 39 & 1369 & 37 & 2 & $K_{39}$ & $(37,28,21)$ & B, c \\
60543 & 5952 & 1147 & 576 & 63 & 961 & 31 & 2 & $(63,32,16,16)$ & $(31,6,1)$ & S \\
60543 & 14880 & 7192 & 3600 & 63 & 961 & 31 & 2 & $(63,32,16,16)$ & $(31,15,7)$ & P,S \\
60543 & 15872 & 8184 & 4096 & 63 & 961 & 31 & 2 & $(63,32,16,16)$ & $(31,16,8)$ & P,S, c \\
60543 & 24800 & 19995 & 10000 & 63 & 961 & 31 & 2 & $(63,32,16,16)$ & $(31,25,20)$ & S, c \\
83205 & 39732 & 19393 & 18963 & 45 & 1849 & 43 & 2 & $K_{45}$ & $(43,21,10)$ & P \\
83205 & 41624 & 21285 & 20812 & 45 & 1849 & 43 & 2 & $K_{45}$ & $(43,22,11)$ & P, c \\
\bottomrule
\end{tabular}
\end{table}